\documentclass{svmult}
\usepackage{amsmath,amssymb}
\usepackage{graphicx}
\usepackage{helvet}       
\usepackage{courier}

\newtheorem{Theorem}{Theorem}
\newtheorem{Proposition}[Theorem]{Proposition}
\newtheorem{Definition}{Definition}

\begin{document}
\title*{The two-dimensional density of Bernoulli Convolutions} 
\author{Christoph Bandt}
\institute{This work was supported by Deutsche Forschungsgemeinschaft, grant Ba 1332/11-1.\\ Christoph Bandt\at Institute of Mathematics, University of Greifswald, 17487 Greifswald, Germany. \\ \email{bandt@uni-greifswald.de} }
\maketitle

\abstract{Bernoulli convolutions form a one-parameter family of self-similar measures on the unit interval. We suggest to study their two-dimensional density which has an intricate combinatorial structure.  Visualizing this structure we discuss
results of Erd\"os, J\'oo, Komornik, Sidorov, de Vries, Jordan, Shmerkin and Solomyak, Feng and Wang. 
We emphasize the r\^ole of finite orbits of associated multivalued maps and mention a few new properties and examples.}
 
\section{Introduction}
The Bernoulli convolution with parameter $\beta\in (1,2] ,$ or $t=\frac{1}{\beta},$ is the unique probability measure $\nu$ on $[0,1]$ which fulfils
\begin{equation} 
\nu (A)=\frac12 \nu(g_0(A))+\frac12 \nu(g_1(A))\quad\mbox{ for all Borel sets } A\subset [0,1] ,\label{bern}\end{equation}
where 
\[ g_0:[0,t]\to [0,1], \ g_0(x)=\beta x\quad\mbox{ and }\quad g_1:[1-t,1]\to [0,1], \ g_1(x)=\beta x+1-\beta\] 
are linear functions with the same slope $\beta$ and fixed points 0 and 1, respectively.  We choose $t=\frac{1}{\beta}$ as parameter in $[\frac12, 1)$ and write  $\nu=\nu_t$ if necessary. If $A$ is not contained in the overlap interval $[1-t,t]$ then either $g_0$ or $g_1$ is not defined on a part of $A,$ and this part is ignored in the corresponding term of \eqref{bern}.
From the viewpoint of fractals, $\nu$ is a self-similar measure with respect to the contractions $f_0(x)=tx, f_1(x)=tx+1-t$ which are the inverse maps of $g_0, g_1.$ The support of $\nu$ is always the unit interval. 

These are the simplest fractal constructions with overlap, and their structure is not yet understood. Only for a countable set of Garsia numbers $\beta$ it is known that $\nu$ has a density \cite{G}. Already 1939 Erd\"os proved that a density of $\nu$ does not exist for Pisot numbers $\beta$ which also form a countable set \cite{E}. (Pisot and Garsia numbers are defined at the end of this section.) For all other parameters $\beta ,$ including all rational numbers, it is not yet known whether $\nu$ is singular or absolutely continuous. However, Solomyak \cite{So} proved 1995 that the set of 'singular parameters' $\beta$ has Lebesgue measure zero. Recently, Shmerkin \cite{Sh14} applied a technique of Hochman to show that this set even has Hausdorff dimension zero. See the surveys \cite{PSS, So4} for more information on the history of Bernoulli convolutions and \cite{Sh14,JSS,FS,F11,HHM} for some recent results. Very recently, Varju \cite{Va} proved that all algebraic numbers $\beta<1+\epsilon ,$ where $\epsilon$ depends on the Mahler measure of $\beta ,$ will lead to absolutely continuous measures. No estimate for $\epsilon$ was given, however. If $\epsilon$ is astronomically small, the result is difficult to interpret since $\nu$ for $\beta\to 1$ converges to the Dirac measure at $\frac12 .$

In this note, we give a non-technical introduction to the combinatorial structure of all Bernoulli convolutions. We focus on computer-generated figures and refer to \cite{BZ} for details. 
Solomyak's theorem, in the $L^2$ version given by Peres and Solomyak \cite{PS}, can be reformulated as follows.

\begin{Theorem}[2D density of Bernoulli convolutions \cite{So,PS}] \label{2D}\hfill\\
There is an $L^2$ function $\Phi: [\frac12,1]\times[0,1]\to [0,\infty)$ such that for Lebesgue almost all parameters $t=1/\beta\in  [\frac12,1],$ the density of the Bernoulli convolution $\nu_t$ is the function $\Phi(t,x), x\in [0,1].$
\end{Theorem}

\begin{figure}[h]
\includegraphics[width=0.999\textwidth]{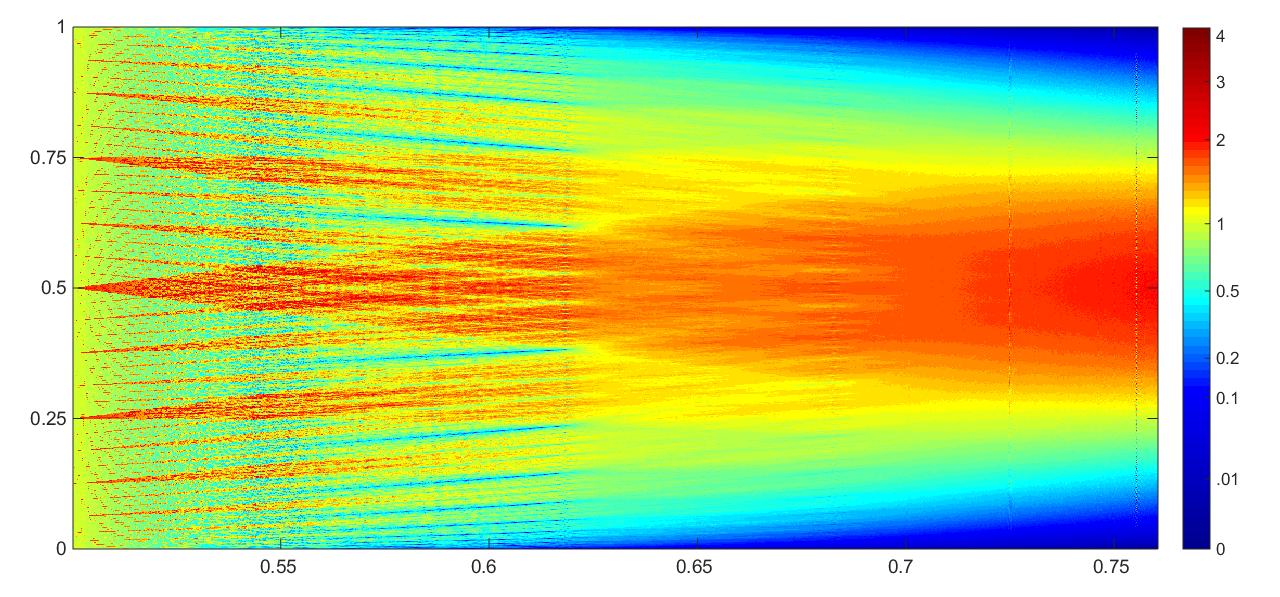}
\caption{The function $\Phi$ for $0.5\le t\le 0.76.$ Bernoulli convolutions for 1000 values $t=1/\beta$
were approximated by histograms with 20000 bins and visualized as vertical sections. The color code on the right indicates that $D=\{ (t,y)|\ t\le y\le 1-t\}$ contains large values. }
\label{B5}
\end{figure}

Thus instead of a bundle of different measures $\nu_t,$ we study one function of two variables describing the whole Bernoulli scenario. For $t$ in $[\frac12 , 0.76],$ the function $\Phi(t,y)$ is sketched in Figure \ref{B5} as color-coded map.
The apparent structure is connected with results of different authors and will be explained below. Our algorithms generating the measures $\nu_t$  include the 'chaos game', inverse iteration \cite[Chapter 8]{bar}, and approximation by Markov chains. 

In Section 2 we start with some carefully calculated histograms of Bernoulli convolutions to get an idea of their properties. The numerical appearance of $\nu_t$ is chaotic when $t$ goes down to $\frac12 ,$ for all parameters, not only Pisot. On the other hand, singular $\nu_t$ for Pisot parameters look more harmless than one would expect.  Obviously there is large variation, and there are also higher peaks, compared with neighboring parameters. However, within ordinary numerical accuracy -- histogram bars of width greater $10^{-8},$ say -- we could not find values larger than 10. 

In Section 3 we explain the basic overlap structure. In Section 4 results of \cite{DaKa,EJK,GS,Si3,Si9,JSS} on points with unique addresses are illustrated with the function $\Phi .$ Section 5 discusses kneading sequences, related to work in \cite{dVK,ACS}. In Section 6, intersections of kneading curves are studied, and results in \cite{FW,F11} improved. For details we refer to \cite{BZ}.

Polynomials of $\beta$ arise as repeated compositions of $g_0,g_1,$ and $\beta$ is a root of a polynomial with integer coefficients if certain higher-level overlaps in the fractal construction coincide. Let us mention some terminology. A root of a polynomial with integer coefficients and leading coefficient one is called an algebraic integer. We consider only positive real roots $\beta .$ There is a minimal polynomial, the other roots of which are called conjugates of $\beta .$ If all conjugates are strictly smaller than one in modulus, $\beta$ is called a Pisot number. If the conjugates' modulus is not greater than one, and equal one for at least one conjugate, $\beta$ is termed Salem number. If the modulus of all conjugates is larger one, and the constant term of the minimal polynomial is 2 or -2, $\beta$ is called a Garsia number.  If $\beta$ is strictly greater than the modulus of all its conjugates, we call $\beta$ a Perron number. When the inequality need not be strict, $\beta$ is a weak Perron number.

\section{Five phases of Bernoulli convolutions}
\begin{figure}[h!] \begin{center}
\includegraphics[width=0.9\textwidth, height=0.24\textwidth]{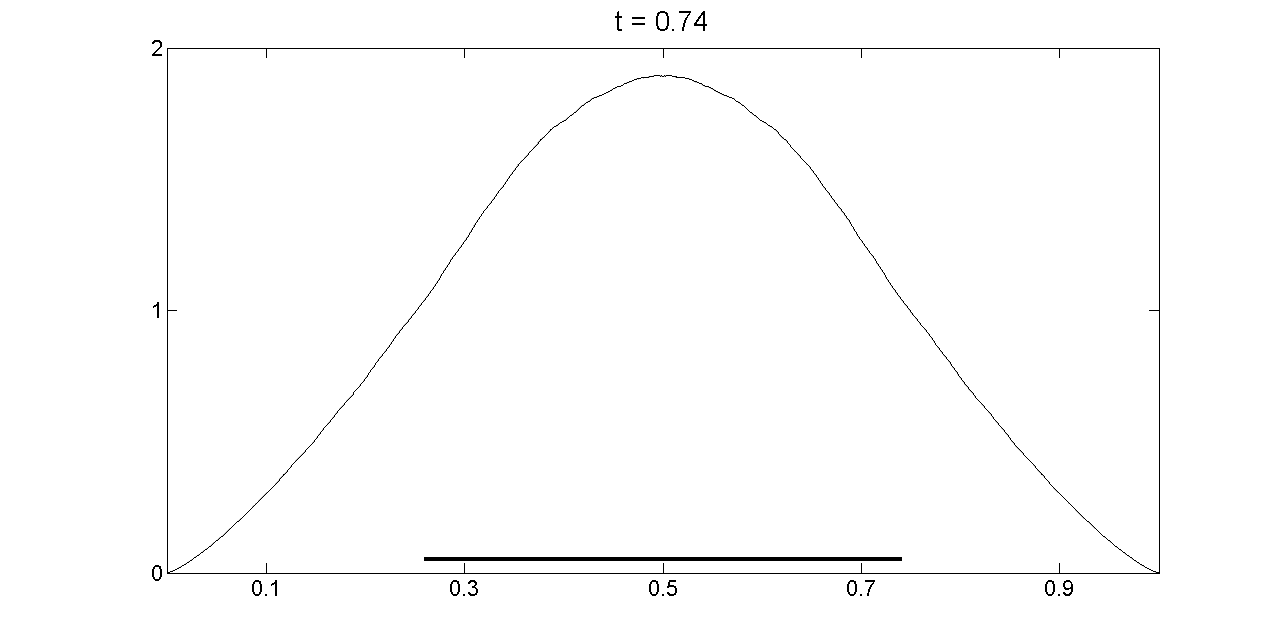}\\
\includegraphics[width=0.9\textwidth, height=0.2\textwidth]{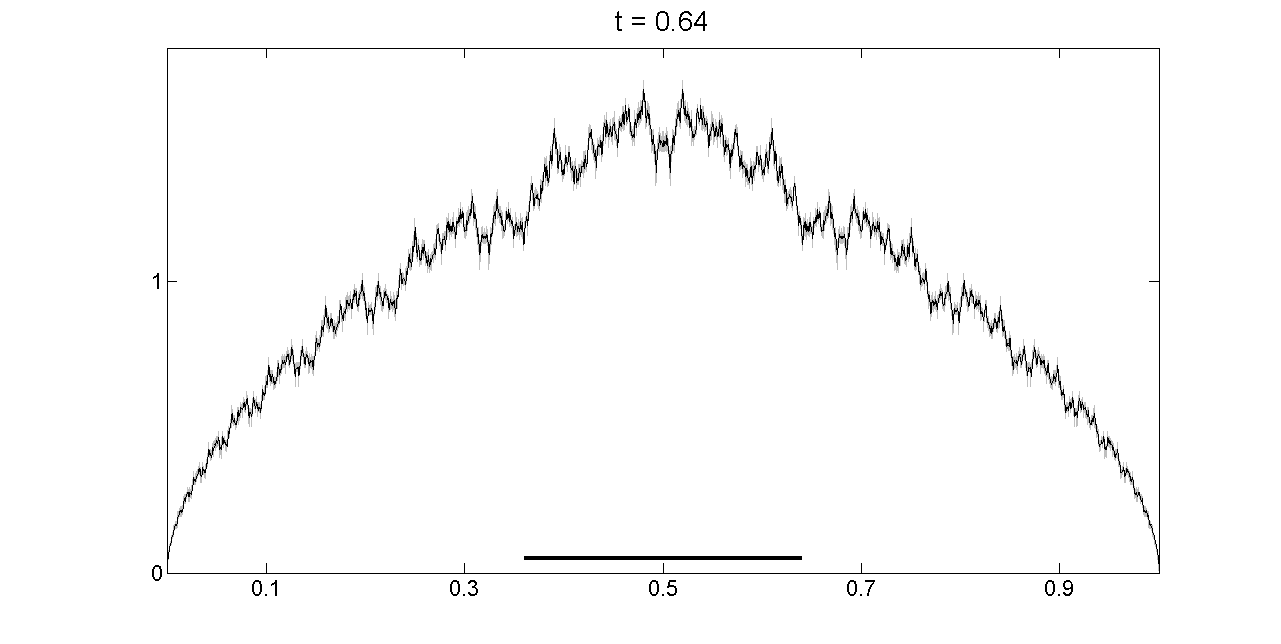}\\
\includegraphics[width=0.9\textwidth, height=0.37\textwidth]{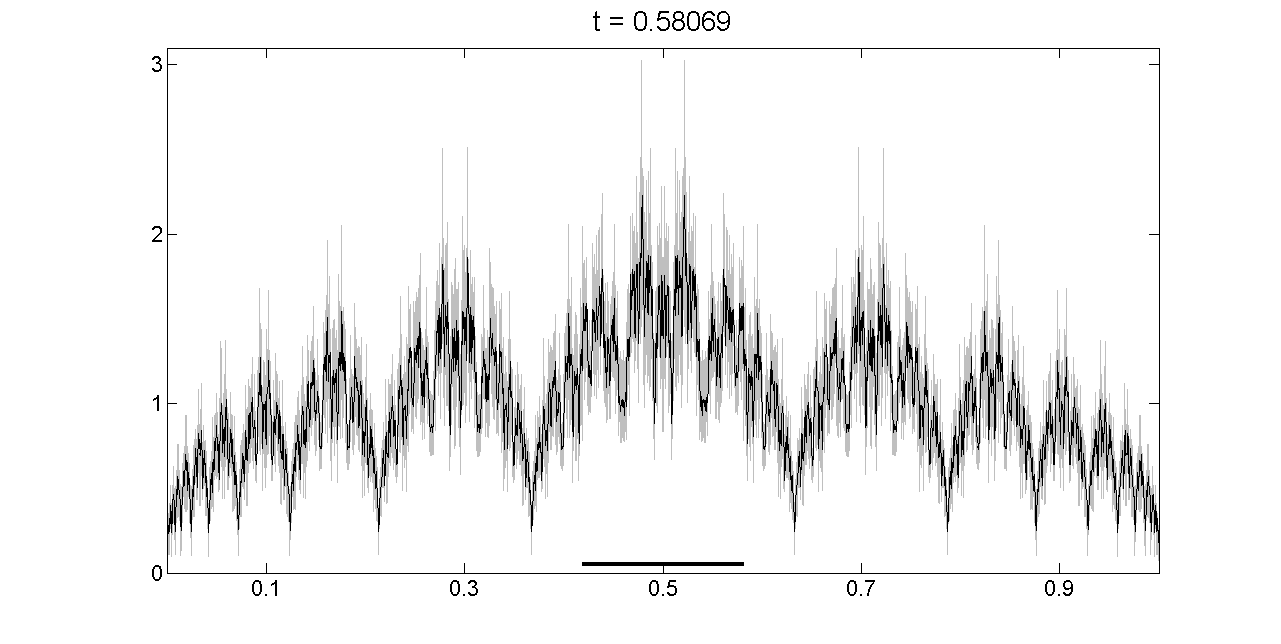} \end{center}
\caption{Bernoulli convolutions for different $t,$ indicating phases 1,2,3.  Histograms with 2000 bins in black and 50000 bins in grey. The bar indicates the overlap interval. The last example is assumed to be singular (Salem number).}
\label{BC1}
\end{figure}

\begin{figure}[h!] \begin{center}
\includegraphics[width=0.9\textwidth, height=0.27\textwidth]{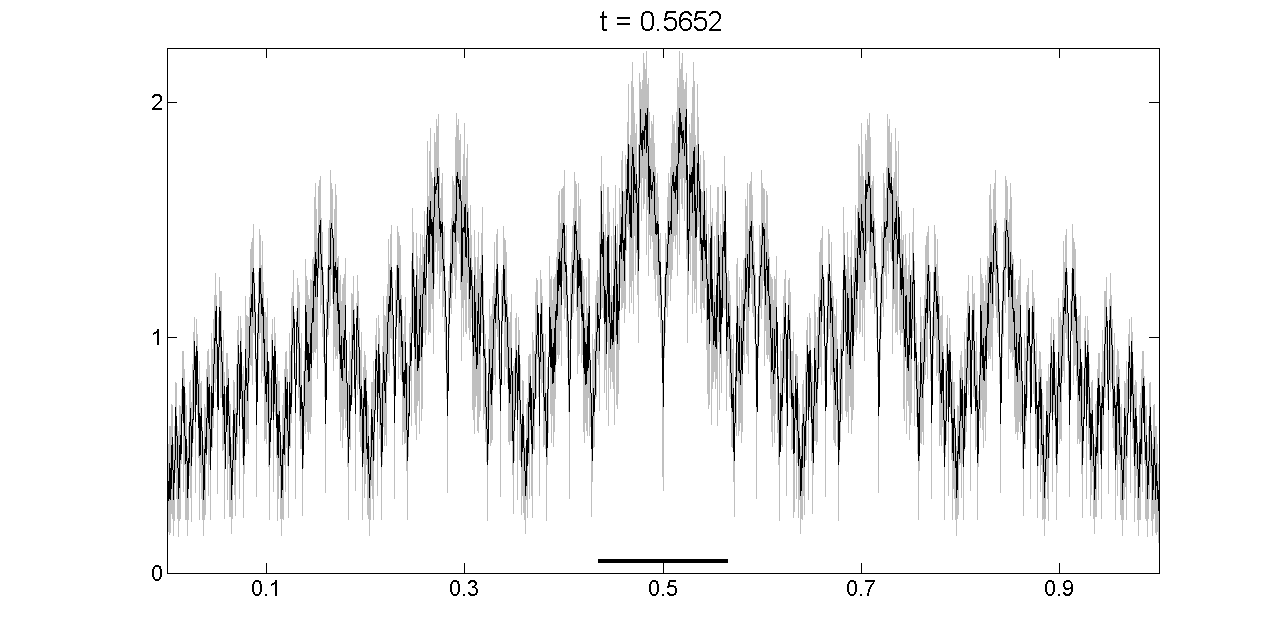}\\
\includegraphics[width=0.9\textwidth, height=0.37\textwidth]{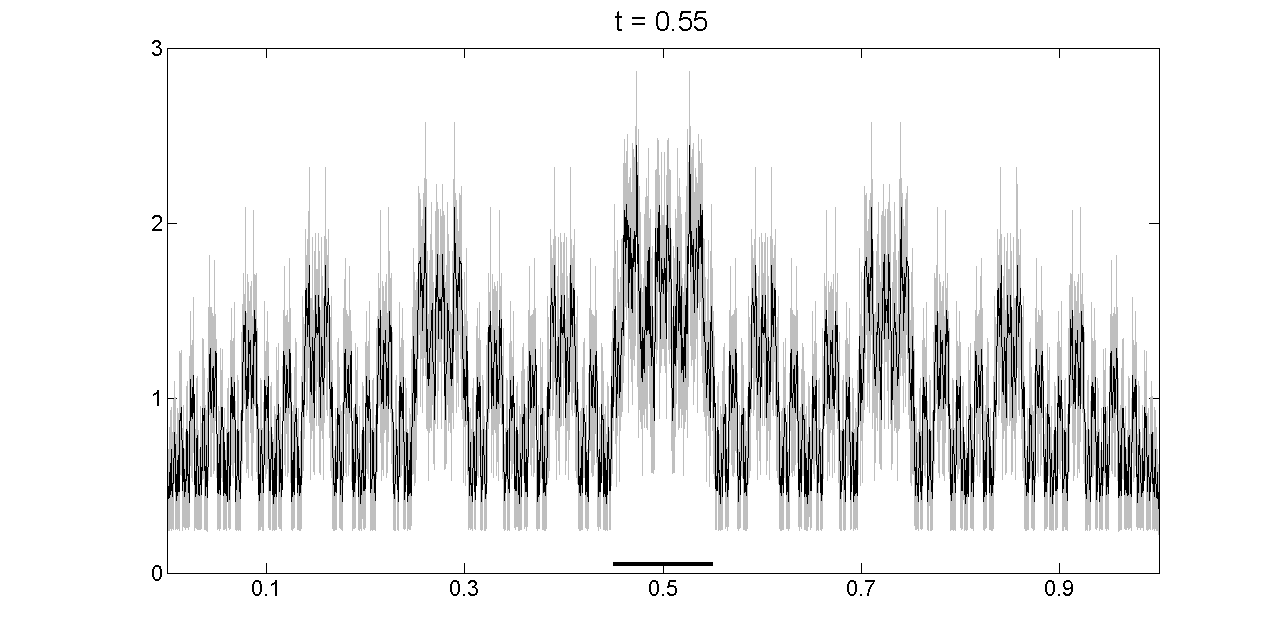}\\
\includegraphics[width=0.9\textwidth, height=0.52\textwidth]{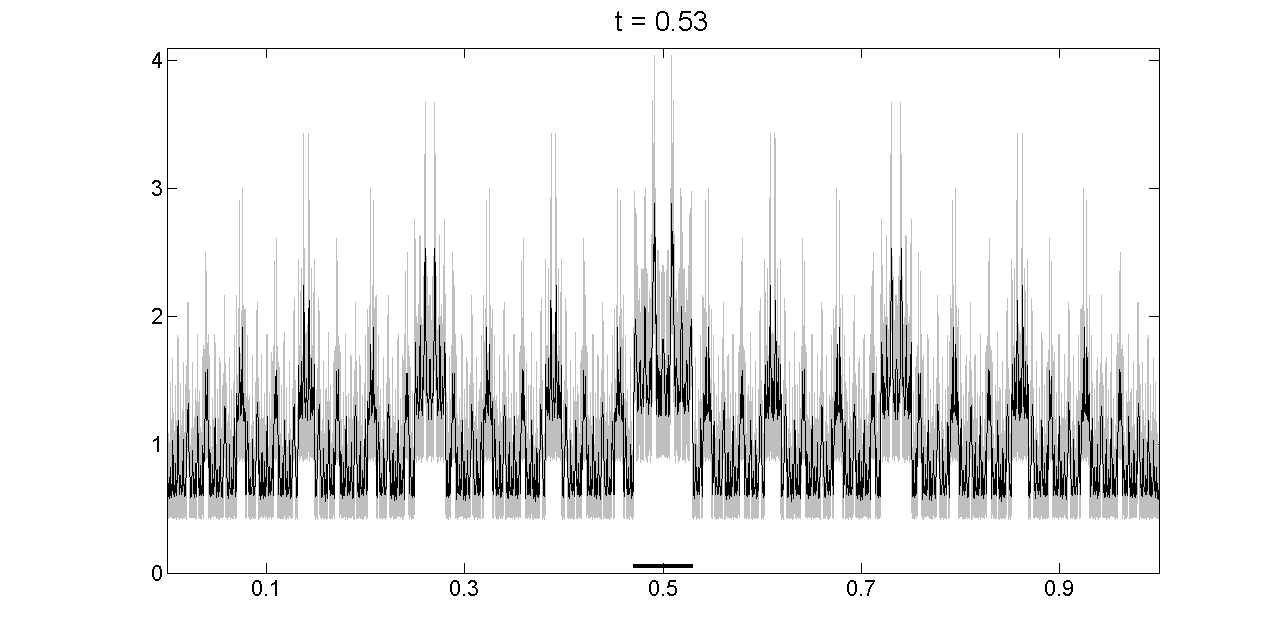} \end{center}
\caption{Bernoulli convolutions for different $t,$ indicating phases 3,4,5.  Histograms with 2000 bins in black and 50000 bins in grey. The bar indicates the overlap interval. The first example is a density function (Garsia number).}
\label{BC2}
\end{figure}

A few pictures of single Bernoulli convolutions will show what kind of vertical sections of $\Phi(t,y)$ are put together. We roughly distinguish five phases.
A specimen for each phase is given in Figures \ref{BC1} and \ref{BC2}. The concept of zero will be made precise below.

\begin{enumerate} \item $1< \beta\le\sqrt{2},$ or $0.707\le t<1.$ The density functions are smooth. Most of them resemble a normal distribution.
\item $\sqrt{2}<\beta<\frac12 (1+\sqrt{5})=\tau ,$  or $0.618<t<0.707.$ The density functions are not smooth, but continuous and strictly positive.
\item $\frac12 (1+\sqrt{5})\le\beta\le\beta_{KL},$ or $0.5595\le t\le 0.618.$ The densities have countably many zeros.
\item $\beta_{KL}<\beta<\tau_3,$ or $0.5437<t<0.5595.$ There are uncountably many zeros outside the overlap interval $[t,1-t]$ and finite or countably many zeros inside that interval.
\item $\tau_3\le \beta<2,$ or $0.5<t\le 0.5437.$ There are Cantor sets of zeros inside and outside the overlap region. The dimension of these Cantor sets approaches one when $t$ goes to 0.5. The density functions seem to have almost vertical slope everywhere.
\end{enumerate}

Our description of phases is not a rigorous mathematical statement.  For certain exceptional numbers $\beta$ it is definitely false. The Fibonacci number $\tau$ and the Tribonacci number $\tau_3=1.8393,$ the root of $x^3-x^2-x-1,$ are Pisot numbers for which no density exists. $\beta_{KL}=1.7872$ is the parameter found by Komornik and Loreti \cite{KL}. Positivity, continuity and smoothness of the functions can be proved only for very particular parameters like $\sqrt{2},$ since in general we do not even know whether a density exists. Sidorov \cite{Si9} called phases 3,4,5 the lower, middle and top order. He proved that in phase 5 there is at least one zero inside the overlap region. In general, there seems to be a Cantor set of zeros. Some other illustrations of Bernoulli convolutions can be found in \cite{BB,So4}.  

Note that zeros of the density functions usually do not exist in a numerical sense. Even if we draw the graph of a density function as a histogram with 10 million bars, there will be no proper zero in phase 5. The assertion on Hausdorff dimension of zeros was shown rigorously in \cite{JSS}. Nevertheless, zeros are so thin that they are not recognized numerically.  

On the other hand, even for parameters where densities cannot be bounded, as on bottom of Figure \ref{BC1}, the maximum values of our functions are between 2 and 4, depending on the resolution of the picture. To study the resolution effect, our histograms were drawn with 2000 bins in black and 50000 bins in grey. 
For phases 1 and 2, differences are hardly visible but they matter: we could not decide numerically for which $t$ near 1 the density is increasing up to $x=\frac12 .$ In phases 3 to 5, resolution differences seem tremendous.  The two examples for phase 3 include a Salem number for which Feng \cite{F11} proved that a bounded density function cannot exist, and a Garsia number with bounded density which equals zero in the central point $\frac12 .$ The Salem parameter leads to  larger peaks and slightly larger variation. 

This discussion shows how difficult it is to study Bernoulli convolutions one by one. The treatment of the two-dimensional density $\Phi$ will be easier.

\begin{figure}[h]
\begin{center}\includegraphics[width=0.9\textwidth]{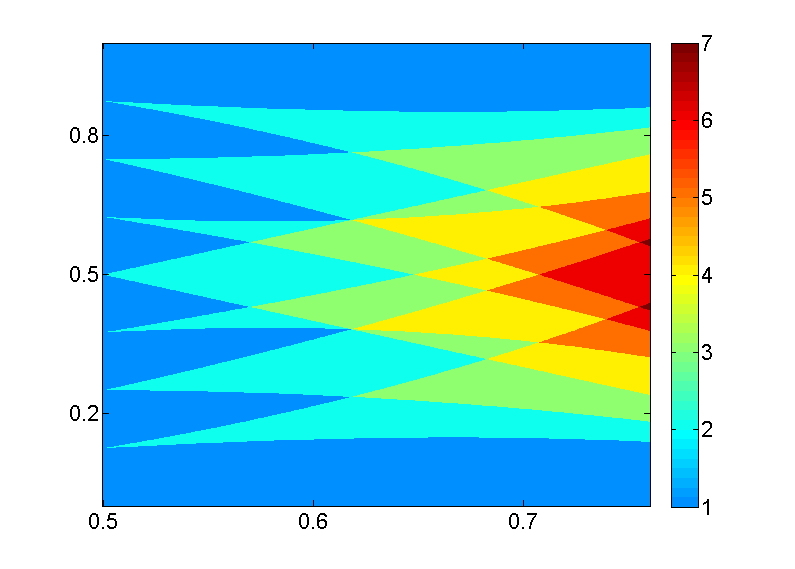}\end{center}
\caption{A low-order approximation of $\Phi .$ There are 0 up to 6 horns of levels 0,1 and 2 which meet in a point $(y,t),$ or 1 up to 7 values of $g_t^3(y)$ in Section 5. Landmark points $\beta$ are algebraic integers like golden mean and $\sqrt{2},$ cf. Table \ref{pigar}. Borders as in Figure \ref{B5}.}
\label{B3}
\end{figure}

\section{Overlap region and horns}
The most obvious feature of Figure \ref{B5} is the overlap region $D=\{ (t,y)|\ t\le y\le 1-t\},$ a big triangle with large values of $\Phi .$ It represents the overlap intervals $D=[1-t,t]$ on the first level of the fractal construction (we should write $D_t$ but omit $t$). In the definition \eqref{bern} of $\nu(A),$ one of the terms on the right-hand side becomes zero if $A\cap D=\emptyset$ so the density is smaller outside the overlap region. However, there are other 'horns' which come out from the big triangle and which represent the overlaps on higher level, for example $D_0=f_0(D)$ which is given by $t^2\le y\le t(1-t).$ The general form is $D_w=f_w(D)$ where $w=w_1...w_n$ is a 0-1-word, and $f_w=f_{w_1}\cdots f_{w_n}.$ 

Since $D_w$ is mapped by $g_w$ on $D,$ the interior structure of the horn $D_w$ reflects the interior structure of $D,$ at least to some extent, according to  \eqref{bern}. The equations of lower and upper border of $D_w$ are $y=f_w(1-t)$ and  $y=f_w(t).$ These are polynomials in $t$ with coefficients $\pm 1$ and zero, already studied by Garsia \cite{G}. Lower borders do not intersect each other and meet in $(1,0).$ Upper borders do not intersect and meet in $(1,1).$ Landmarks are obtained from intersection points of lower and upper borders of different horns. The corresponding parameters $\beta$ are algebraic integers.

Figure \ref{B3} shows the seven horns $D, D_0, D_1, D_{00}, D_{01}, D_{10}, D_{11}.$ Up to six horns intersect in a point $(t,y)$ with $t\le 0.76$ which means $y$ has up to 7 addresses when only three levels of iteration are studied. 
The two parameters of degree two, obtained from the intersection of $D_0$ with $D$ and $D_1,$ are the golden mean at $t=0.618,$ and the Garsia number $\sqrt{2}, t=0.707.$ 
Further landmarks $\beta$ in Figure \ref{B3} are on the curve $y(t)=t-t^2+t^3=f_0f_1(1-t)$ describing the upper border of $D_{01}$ with tip at $y(\frac12)=\frac{3}{8}.$
The intersection with curves $1-t$ and $1-t^2$ of lower order horns $D, D_1$ leads to well-known Pisot numbers at $t=.570$ and $t=.6823$ while intersections with horns $D_{10}, D_{11}$ of the same order yields Garsia numbers at $t=.648$ and $t=.739.$ Parameters in the last row of Table \ref{pigar} come from intersections of $D_{00}$
with $D_{1}$ and $D_{11}.$

\begin{table}[h!] \begin{center}\normalsize
\begin{tabular}{|cc|cc|} \hline
\multicolumn{2}{|c|}{Pisot numbers}& \multicolumn{2}{c|}{Garsia numbers}\\ \hline
$t=1/\beta$&$p(\beta)$& $t=1/\beta$&$p(\beta)$\\ \hline \hline
 .618&$\beta^2-\beta-1$& .707&$\beta^2-2$\\ \hline
 .570&$\beta^3-2\beta^2+\beta-1$& .648&$\beta^3-2\beta^2+2\beta-2$\\ \hline
 .682&$\beta^3-\beta^2-1$& .739&$\beta^3-\beta^2+\beta-2$\\ \hline
 .755&$\beta^3-\beta-1$& .794&$\beta^3-2$\\ \hline
\end{tabular} \end{center}
\caption{Landmarks in Figure \ref{B3}. Pisot parameters include golden mean, its doubling counterpart, plastic number, first Pisot number. Garsia polynomials obtained by adding $\beta -1.$}\label{pigar}
\end{table}

On this level, all landmark points correspond to the two classes of numbers which have been thoroughly studied in connection with Bernoulli convolutions. Moreover, polynomials of Garsia numbers are obtained by adding $\beta -1$ to the corresponding Pisot polynomial. On higher level, the situation is more complicated. Not all horns will intersect, and we often get only Perron numbers \cite{BZ}.

In the sequel we focus our study on two regions: those points which are not contained in any horn, and the points inside $D.$ Due to symmetry, it suffices to consider points $(t,y)$ with $y\le\frac12 .$ In this note, we shall concentrate on phases 3 to 5 where the structure of $\Phi$ is most apparent. Figure \ref{B6} shows a magnification of this part of Figure \ref{B5}.

\section{Points with unique addresses}
Already before 1990, it was discovered that for $t<0.618$ there are points which have a unique address in the fractal construction. For $t>0.618$ all points $y\in [0,1]$ have a continuum of addresses while for the golden mean parameter $t=0.618$ all points have infinitely many addresses, where a countable dense set of points has only a countable number of addresses. See Daroczy and Katai \cite{DaKa}, Erd\"os, Jo\'o and Komornik \cite{EJK}, Glendinning and Sidorov \cite{GS}, Sidorov \cite{Si3,Si9} and various other papers quoted there.

An address of $y$ is a 01-sequence $s=s_1s_2...$ with $y=\lim_{n\to\infty} f_{s_1}\cdots f_{s_n}(x_0)$ where the initial point $x_0$ does not matter.  Points in the overlap region have at least two addresses, and so do all points in all horns $f_w(D).$ \emph{Thus the points with a unique address coincide with the points which do not belong to any horn.}

\begin{figure}[h]
\includegraphics[width=1.03\textwidth]{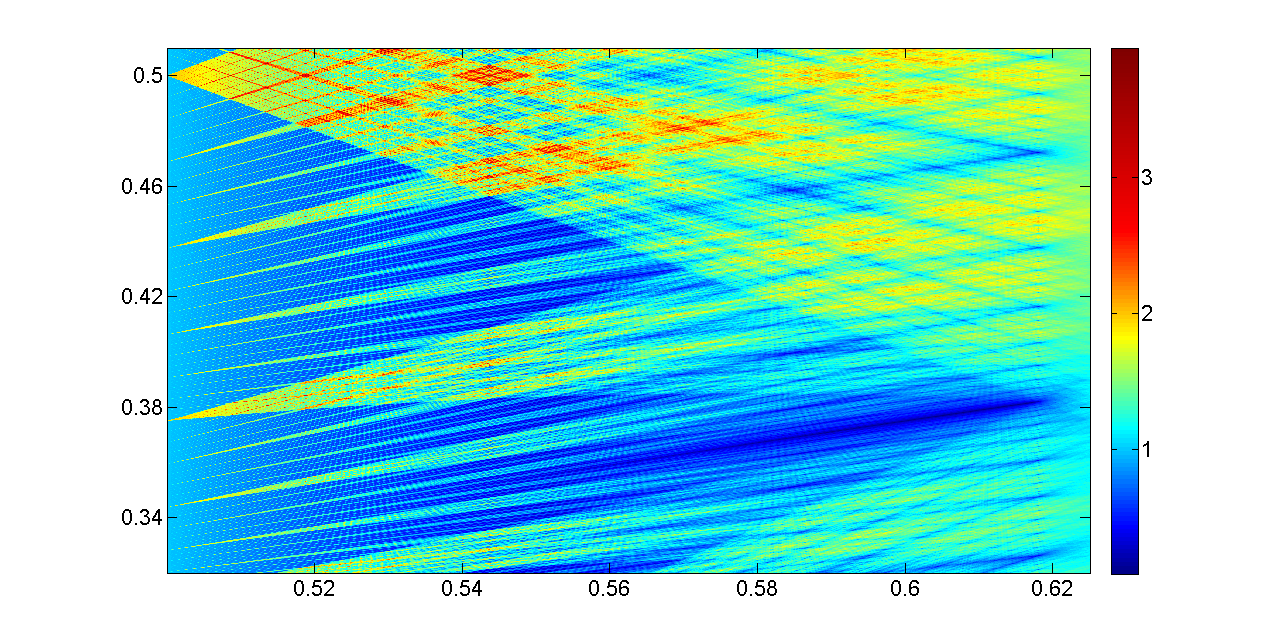}
\caption{$\Phi$ for $0.5\le t\le 0.63,\quad 0.325\le y\le 0.5.$ Right of the golden mean, the structure becomes blurred. At $t=0.618,$ a periodic point appears outside $D,$ describing a dark curve which hits the $y$-axis at $1/3.$}
\label{B6}
\end{figure}

In Figure \ref{B6} these points are recognized by their dark color. Why? Because the measure of an $\epsilon$-neighborhood $B(y,\epsilon)=[y-\epsilon, y+\epsilon]$ of a point $y$ with unique address decreases very fast with $\epsilon .$ The appropriate parameter is \emph{local dimension.} Roughly speaking, a measure $\nu$ has dimension $d=d_y(\nu )$ at a point $y$ if $\nu (B(y,\epsilon))\approx\epsilon^d$ for small $\epsilon ,$ and Lebesgue measure on $\mathbb{R}^d$ is a basic model for this concept. The precise definition is 
\begin{equation} d_y(\nu )= \lim_{\epsilon\to 0} \frac{\log \nu (B(y,\epsilon))}{\log \epsilon} \ .
\label{loc}\end{equation}
We consider only cases where the limit exists, so we need not distinguish upper and lower local dimension. See \cite{fal} for details.

If for a Bernoulli measure $y$ has only one address, then on each level $n$ of the fractal construction, $y$ is contained in a single piece, an interval of  length $t^n$ and $\nu$-measure $2^{-n}.$ Putting this interval as approximation of $B(y,t^n)$ into \eqref{loc} we obtain
$d_y(\nu_t)=\frac{\log 2}{\log \beta}>1 .$ For proofs of related statements, see \cite[Proposition 4]{BZ}, \cite{FS}. 

A more basic concept is the density of $\nu$ at $y,$ 
\[ D_\nu(y)=\lim_{\epsilon\to 0} \frac{\nu (B(y,\epsilon))}{2 \epsilon}\ .\]
This function of $y$ is called the density function of $\nu$ when the limit exists for all $y.$ 
If the local dimension of a measure $\nu$ at $y$ is greater than 1, then clearly $D_\nu(y)=0 .$ \smallskip

\emph{Conjecture. } 
For Bernoulli measures $\nu_t$ the converse seems to be true: density zero implies local dimension greater than 1.
We further conjecture that with exception of weak Perron parameters $\beta$ discussed in  \cite[Theorem 6]{BZ},  the \emph{only} value of local dimension larger one of a Bernoulli convolution equals $\frac{\log 2}{\log \beta}.$ \smallskip

This value describes all cases where $y$ has finite or countable number of addresses. Certainly $\frac{\log 2}{\log \beta}$ is the largest possible local dimension since every point has at least one address. Sidorov and Baker studied various examples of points with two, three or countably many addresses \cite{Si9,Baker,BS,BZ}, but points with unique address are much more frequent.

Local dimension explains why points with unique addresses are so apparent in the graph of $\Phi.$  Jordan, Shmerkin and Solomyak \cite[Theorem 1.5]{JSS} found sets $\tilde{A}_t$ of local dimension $\frac{\log 2}{\log \beta}$ for $\frac12<t<0.618,$ and proved that their Hausdorff dimension tends to 1 for $t\to\frac12 .$ As it turns out, the  $\tilde{A}_t$ contain exactly the points with unique addresses studied by other authors.\smallskip 

Figure \ref{B6} shows that the value $\frac{\log 2}{\log \beta}$ matters a lot. At $t=0.618$ this value is 1.44 and there is a very thick dark line in Figure \ref{B6} which corresponds to a broad valley of the density, as on the bottom of Figure \ref{BC1}. For $t\approx\frac12$ the largest local dimension approaches 1. This implies average coloring in Figure \ref{B6} and extremely narrow valleys on the bottom of Figure \ref{BC2}. The points 0 and 1 have a unique address for every $t,$ which yields the dark margin of Figure \ref{B5}. For $t\to 1$ the local dimension at 0 and 1 converges to infinity. \smallskip

\emph{Remark. } At all points $(t,y)$ with unique address, the local dimension of the two-variable function $\Phi$ equals $1+\frac{\log 2}{\log \beta}.$  Thus $d_{(t,y)}(\Phi)$ assumes all values between 2 and infinity.  \smallskip

\emph{Sketch of proof. } For the second assertion, it is enough to consider $y=0$ for all $t.$ So let us prove the first assertion only for this case.
Local dimension of $\Phi$ is defined with the associated measure $\mu (A\times B)=\int_A \int_B \Phi (t,y) d\nu_t(y) dt$.
Take a point $(t,0)$ with $\frac12 <t<1.$ Choose $\epsilon$ so that $A=[ t-\epsilon, t+\epsilon]$ is a subset of $(\frac12 ,1)$ and $B= [0,\epsilon]$ fulfils
$\frac12 \nu_t(B) \le  \nu_s (B)\le 2 \nu_t (B)$ for all $s\in A .$ Then \
$\epsilon \nu_t(B)\le \mu(A\times B)\le 4 \epsilon \nu_t(B)$ \  and 
\[
\lim_{\epsilon\to 0} \frac{\log \mu(A\times B)}{\log\epsilon} =1+ \lim_{\epsilon\to 0} \frac{\log \nu_t (B)}{\log \epsilon}
=1+\frac{\log 2}{\log \beta}\ . \]
 \smallskip

For a general proof we need the fact that each $(t,y)$ with $t>\frac12$ and unique address lies on a differentiable curve $y_b(t).$ This will be shown below.

\section{Quantile curves and conjugacy with the doubling map}
Here we stress the fact that \emph{points with unique addresses appear on differentiable curves and thus constitute a smooth element in an otherwise chaotic scenario.} We start with an example.

The remarkable phase transition at 0.618 is essentially caused by a single periodic orbit $x_0=\frac{t}{1+t},\,
x_1=\frac{1}{1+t}.$  We have $g_0(x_0)=x_1$ and $g_1(x_1)=x_0$ for every $t.$ Thus if $x_0,x_1$ are outside the overlap triangle -- and this happens exactly for $t<.618$ -- then $x_0$ has the unique address $0101...=\overline{01}$ and $x_1$ has unique address $\overline{10}.$  Thus the smooth function $y(t)=\frac{t}{1+t}$ consists for $\frac12\le t<.618$  entirely of points with unique addresses and provides the prominent curve seen in Figure \ref{B6}.

From the definition \eqref{bern} and the period 2 property of $x_0,x_1$ it follows that $\nu_t [0,x_0]=\nu_t[x_0,x_1]=\nu_t[x_1,1]=\frac13 ,$ cf. \cite[Proposition 2]{BZ}. Thus, in probabilistic terms, the function $y(t)=\frac{t}{1+t}$ defines the $\frac13$-quantile of all measures $\nu_t$ with  $\frac12\le t<0.618.$ This fact will now be generalized.

\begin{Definition}[Doubling map and Bernoulli map]\hfill\\ 
The function $g(x)=2x\mod 1$ defined for $0\le x\le 1$ is called the doubling map. For each $t\in(\frac12, 1)$ the function $g_t(x)=\{ g_0(x), g_1(x)\}$ is the Bernoulli map. It is two-valued on $D=[1-t,t].$ The map $F_t(x)=\nu_t [0,x]$ is the cumulative distribution function of $\nu_t.$  
\end{Definition}

\begin{Proposition}[Conjugacy of Bernoulli and doubling map]\hfill\\ 
For each $t\in (\frac12,1)$ the function $F_t$ defines a conjugacy between the action of $g_t$ on $[0,1]\setminus D$ and the doubling map $g$ on a corresponding subset of $[0,1].$ That is, 
\[ F_t\cdot g_t(x)=g\cdot F_t(x)\qquad\mbox{ for }\quad x\in [0,1]\setminus D\, .\]
\end{Proposition}

This is proved by applying \eqref{bern} to $B=[0,x],$ see \cite{BZ}. To get a conjugacy between dynamical systems, we have to restrict ourselves to points $x$ for which the orbit under $g_t$ does not intersect $D.$ This means $x$ has unique address, thus $t\le 0.618.$ For the borderline cases $x=1-t$ and $x=t,$ the images $g_1(1-t)=0$ and $g_0(t)=1$ will be neglected, to get a closed domain for $g_t.$

\begin{Definition}[Binary itineraries and kneading sequences]\hfill\\ 
A 01-sequence $b_1b_2...$ and the corresponding binary number $b=.b_1b_2...$ are called kneading sequence with respect to the doubling map if no number  $g^{(k)}(b)=.b_kb_{k+1}...$ with $k=1,2,...$ is nearer  to $\frac12$ than $b.$ All preimages $x\in g^{-m}(b)$ of a kneading sequence $b$ for some $m=1,2,...$ are called itineraries. The function $y_b(t)=\frac{1-t}{t}\cdot \sum_{k=1}^\infty b_kt^k$ is called the address curve corresponding to $b.$
\end{Definition}

Itineraries and kneading sequences were introduced in the context of one-dimensional dynamics by Milnor and Thurston
\cite{MT} in slightly different form. Binary
itineraries are exactly those 01-sequences which do not contain $n$ consecutive equal symbols 0 or 1, for some $n.$ The corresponding kneading sequence is obtained by determining the orbit closure of $.b_1b_2...$ under the doubling map, and taking the point (or one of the two points) nearest to $\frac12 .$ If $b$ is a kneading sequence, then so is $1-b.$ So it suffices to study $b<\frac12$ or $b_1=0.$

Similar functions were introduced by Milnor and Thurston to determine the topological entropy of unimodal maps.
The standardizing factor $\frac{1-t}{t}$ comes from our choice of mappings $g_0,g_1$ which define all measures on $[0,1].$ The following theorem says that all address curves with $b_1=0$ can be seen in Figure \ref{B5} as parallel blue curves.

\begin{Theorem}[Address curves define quantiles \cite{BZ}]\hfill\\ 
For each itinerary $b=.b_1b_2...,$ the address curve $y=y_b(t)$ describes the points for which $F_t(y)=b,$ for all Bernoulli measures $\nu_t$ with $\frac12\le t\le t^*$ where $t^*>\frac12$ depends on $b.$
\end{Theorem}

If $b$ itself is a kneading sequence with $b_1=0,$ then $t^*$ is simply the solution of $y_b(t)=1-t.$ For kneading sequences $b$ with $b_1=1,$ we have to solve the equation $y_b(t)=t.$ If $b$ is not a kneading sequence, we determine $t^*$ from the corresponding kneading sequence. In all cases, $t^*$ marks the right endpoint of the dark curve $y_b(t)$ in Figure \ref{B5}.

In the setting of one-dimensional unimodal maps, itineraries are the addresses of points, and kneading sequences are the addresses of the critical point. A 01-sequence is an addresss for different points in different maps, but in a well-behaved parametric family of maps it will appear only once as a kneading sequence, and then it disappears. For Bernoulli convolutions we have a similar situation. Itineraries are unique addresses of certain points, describing the quantile given by the binary number $b.$ At the point $t^*$ they become critical which means that the corresponding kneading sequence is a boundary point of the overlap interval $D=[1-t,t].$ When they enter $D,$ they cease to have a unique address, and to be points of minimal local dimension of $\nu_t.$ As we shall see, the curves $y_b(t)$ remain important beyond $t^*$ although in Figure \ref{B5} $y_b(t^*)$ seems to be their endpoint.

Early results on connections between kneading sequences and Bernoulli convolutions, mostly formulated in the setting of $\beta$-expansions, include the calculation of the Komornik-Loreti parameter \cite{KL} which corresponds to the Feigenbaum point, the description of parameters with unique $\beta$-expansion by de Vries and Komornik \cite{dVK}, and the detection of the Sharkovskii ordering in the set of periodic unique $\beta$-expansions by Allouche, Clarke and Sidorov \cite{ACS}. 

There is a clear one-to-one correspondence between itineraries of unimodal maps, and the unique addresses considered here. To each quadratic map $q_r(x)=rx(1-x)$ without stable periodic points there is a unique corresponding Bernoulli measure $\nu_{t^*}.$ 
The periodic windows of the 'Feigenbaum diagram' \cite{CE} correspond to the horns of the two-dimensional Bernoulli density which intersect the central horn $D.$ The dynamic phenomena are quite different, however. Periodic windows correspond to stable periods while the expanding maps $g_0,g_1$ do not admit stable periodic points. As a consequence the width of horns near the Komornik-Loreti parameter decreases like a double exponential, while the length of periodic windows near the Feigenbaum point follow a famous asymptotically geometric sequence \cite{CE}. This particular detail is hardly visible in our figures. Many other details, as for instance windows of higher order, are more apparent here. When landmark points are determined for the quadratic family, we have polynomials of degree $2^k$ for the $k$-fold iteration of $q_r$ while we get polynomials of degree $k$ for the iteration of $g_t.$

\begin{figure}[th]
\includegraphics[width=0.9\textwidth]{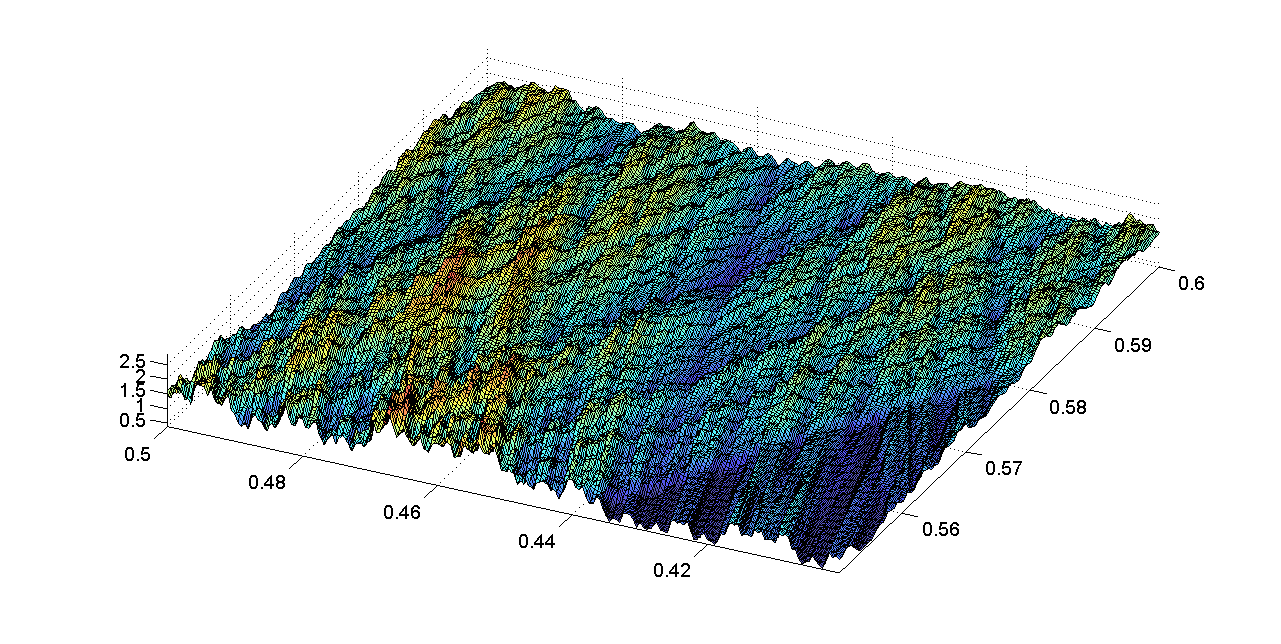}\\
\includegraphics[width=0.9\textwidth]{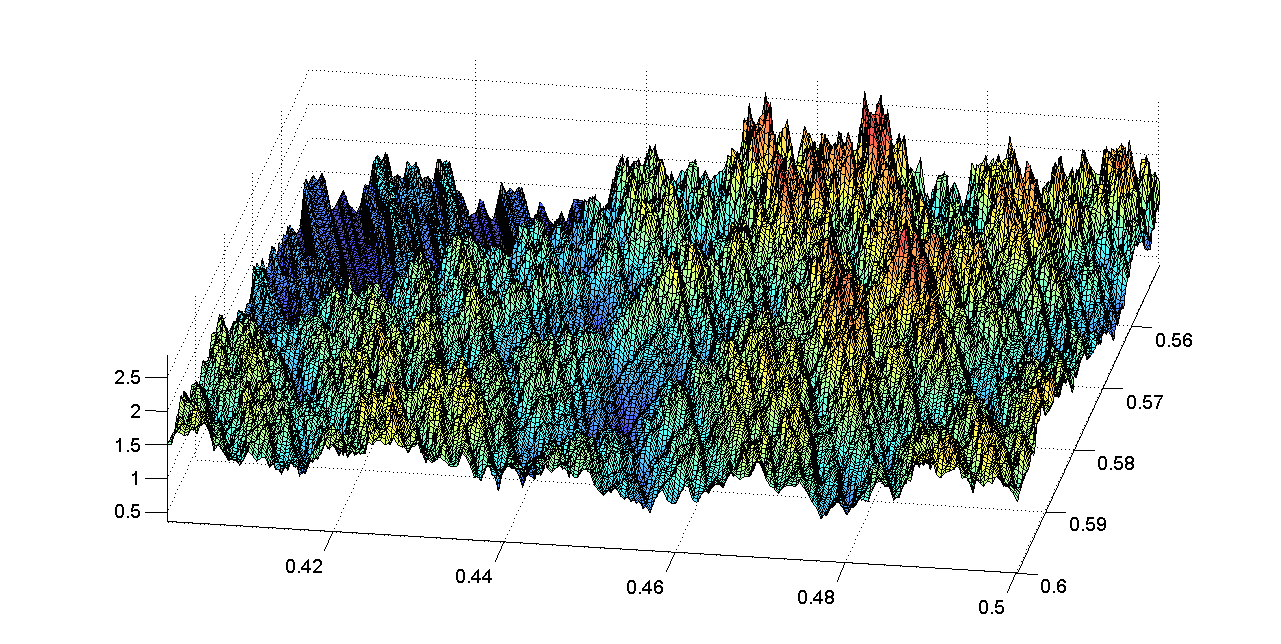}
\caption{Three-dimensional representation of $\Phi$ for $0.55\le t\le 0.6,\quad 0.4\le y\le 0.5$ from two different viewpoints, cf. Figure \ref{B6}. The quantile curves on both sides of the overlap region $D$ act like two wavefronts which interfere in $D.$  There are many valleys and few peaks.}
\label{B7}
\end{figure}

This correspondence seems to deserve more attention. The question for the existence of an absolutely continuous invariant measure of $q_r,$ for example, is wide open as the question for density functions of Bernoulli measures. There could be some connection. For instance, there are parameters with known absolutely continuous measure for $q_r$ which correspond to Garsia numbers. One example is given below. It should also be mentioned that Tiozzo \cite{Ti} recently related the combinatorial structure of kneading sequences to continued fractions and the Gauss map.

Itineraries correspond only to a nowhere dense set of parameters $r$ of quadratic maps, and $t^*$ of Bernoulli measures. One difference between the Bernoulli scenario and real quadratic maps is that there are many dynamic changes between parameters $t^*$ inside the horns, as we shall see below, while there are no dynamic changes within periodic windows of the quadratic family.  In this respect the Bernoulli scenario rather resembles the Mandelbrot set and its abstract model, Thurston's quadratic minor lamination, where windows are replaced by bubbles with rich boundary structure. The overlaps in the Bernoulli scenario lead to other phenomena, however.

\section{Inside the overlap region}
So far we have dealt with the simple part of the Bernoulli scenario. Inside the horn $D$ and its copies $D_w$ we see the chaotic part.  Figure \ref{B7} indicates that the curves $y_b(t)$ remain structure-forming elements of $\Phi .$ However, there are two families of curves which meet inside the horn $D :$ the increasing functions $y_b(t)$ with $b<\frac12 ,$ seen in Figure \ref{B5}, and the family of decreasing functions $y_c(t)$ with 
$c=1-b>\frac12 .$  They look like two wavefronts which interfere with each other. Actually, the situation is more complicated since two families can already meet inside the smaller horns $D_w,$ and then these horns will meet $D,$ and meet each other within $D.$  (See the upper picture of Figure \ref{B7}. The wavefront in the foreground which approaches the line $y=1-t$ does not consist of parallel waves. In the picture below these waves are in the background.)
As a result the continuity of curves is lost. There are fractal mountains. We shall give an explanation for the impression that there are rather few and rather isolated peaks or clusters of peaks while there are many valleys.

A rigorous study of the intersections of itinerary curves confirms this observation. To obtain large values of $\Phi$ in the intersection point, rather restrictive conditions must be fulfilled. We must have two different periodic addresses, which implies that the corresponding parameter $\beta$ is a weak Perron number, and the growth of addresses must be sufficiently fast. We formulate the statement in non-technical form, refer to \cite{BZ} for more details, and give a few examples to clarify the situation.

\begin{Theorem}[Intersections of kneading curves]\hfill\\ 
Let $b=.0b2b3...$ and $c=.1c_2c_3...$ be itineraries, and let the two curves $y_b(t), y_c(t)$ intersect in the point $(s,z)$ inside the central horn $D.$ For (i) and (ii) we assume that no point of the forward image of $z$ under $g_s$ lies in $D.$
\begin{enumerate}
\item[(i)] 
If both $b$ and $c$ have infinite or preperiodic orbit with respect to the doubling map, then $z$ has two addresses, and the local dimension of $\nu_s$ assumes the maximum value $\frac{\log 2}{\log 1/s}$ (except when two 'pre-periods' coincide). 
\item[(ii)]
If one sequence is periodic and the other one is preperiodic or infinite, then $z$ has a countable number of addresses,
and the local dimension is $d_z(\nu_s)=\frac{\log 2}{\log 1/s}.$
\item[(iii)]
If both $b$ and $c$ are periodic with respect to the doubling map, with periods $m$ and $n,$ then $z$ has  uncountably many addresses. If $(2s)^{-m}+(2s)^{-n}>1,$ then the
local dimension $d_z(\nu_s)$ is smaller than 1, and $\nu_s$ cannot have a bounded density.
\end{enumerate}\label{cross}
\end{Theorem}

We illustrate the statement with periodic and preperiodic kneading sequences. Thus the orbit of $z$ under $g_0,g_1$ will be finite. In this case $1/s$ must be a weak Perron number \cite{BZ}. 
The assumption that the orbit remains outside $D$ is restrictive. For $.570<t<.618,$ for instance, we have only the kneading sequences $\frac13 =.\overline{01}$ and $\frac23 =.\overline{10}$  (that is, $t^*<.57$ for other kneading sequences). Thus preperiodic itineraries have the form $,w\overline{10}$ with a 01-word $w,$ corresponding to the binary representation of $k/(3\cdot 2^n)$ for some integer $k.$\smallskip

{\it Examples for case (i). } The curve $y_b(t)=t-t^2+\frac{t^2}{1+t}$ for $\frac{5}{12}=b=.01\overline{10}$ can be seen inside $D$ in Figure \ref{B6} as a dark broken line. The reason is the above theorem: there are so many crossing points with other curves  coming from above which force the density of $\nu_s$ to be zero.  One of these curves, $y_c(t)$ with $\frac{13}{24}=c=.100\overline{01}$ intersects $y_b(t)$ at $t=.618, z=.472.$ This curve can also be followed inside $D$ because of its zeros. The curve $y_d(t)$ with $\frac{25}{48}=c=.1000\overline{01}$ is less apparent but it has a very clear intersection point with $y_b(t)$ at $s=.585, z=.459$ visible in Figure \ref{B6} as a dark crossing. Sidorov \cite{Si9,BS} found this example and proved that $\beta=1/s$ is the smallest parameter where a point has exactly two addresses.\smallskip 

{\it Example for case (ii). }
There are many other points on the curve $y_b(t)$ with $s<t<.618$ where crossings can be seen. They all correspond to case (ii) in Theorem \ref{cross} where we have a countable number of addresses. To give just one example,  
$\frac{8}{15}= c=.\overline{1000}$ with $y_c(t)=(1-t)/(1-t^4)$ leads to $s=.592, z=.463.$\smallskip

{\it Bernoulli convolution with zero at $y=\frac12$. }
We combine $\frac{11}{24}=b=.011\overline{10}$ with $\frac{13}{24}=c=.100\overline{01}.$
When $c=1-b,$ we obtain $s$ as the root of $y_c(t)=\frac12=z .$ In our case we get $s=.565,$ and $\beta=\frac1s$ is a Garsia number with minimal polynomial $\beta^3-2\beta -2.$ This is the smallest $\beta$ for which the central point $\frac12$ has only two addresses. Zero density is visible in the top picture of Figure \ref{BC2}, and the dark crossing is  apparent in Figure \ref{B3}. At first glance, one would think that $\frac12$ usually is a point of maximal density of $\nu_t.$ Except for phase 1, however, this case is indeed an exception.

\begin{figure}[h]
\includegraphics[width=0.999\textwidth]{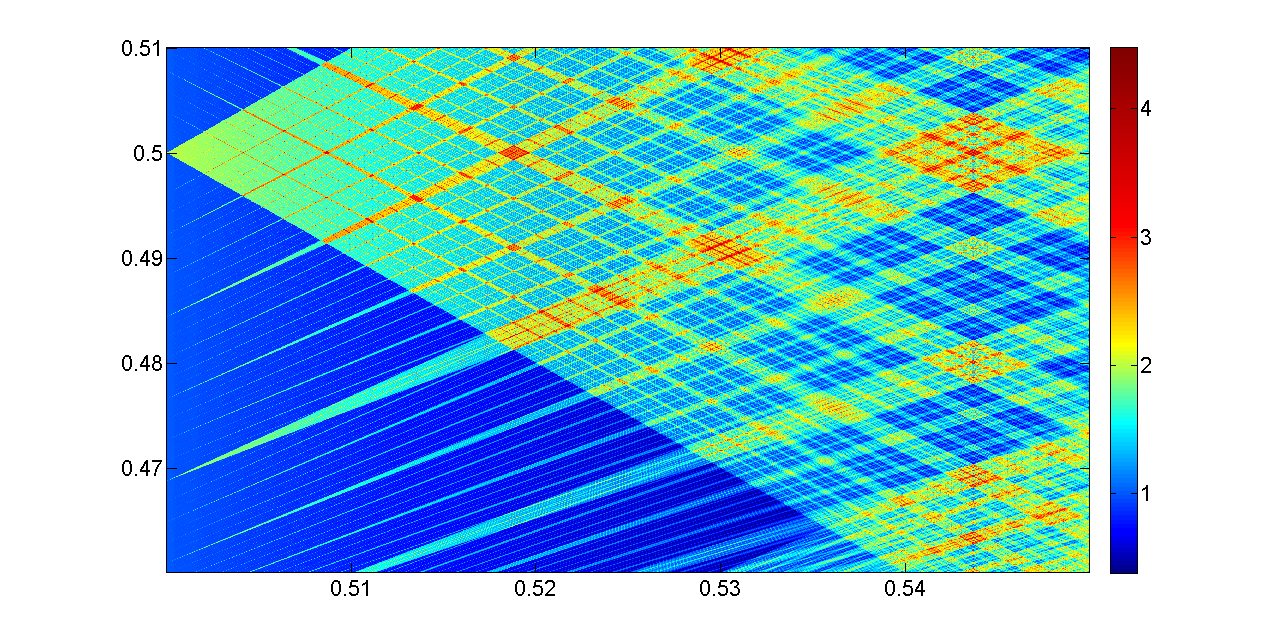}
\caption{Phase 5: $\Phi$ for $0.5\le t\le 0.55,\quad 0.46\le y\le 0.51.$  Even near the Tribonacci parameter on the right, the 
pattern of peaks seems discrete. At the tip $t=y=\frac12$ of the overlap triangle $D$ the scenario looks self-similar, with scaling factor two.}
\label{B8}
\end{figure}

Figure \ref{B8} shows that in phase 5 dark crossings are abundant and form Cantor carpets instead of isolated spots. The reason is the abundance of infinite orbits outside $D.$ All these examples belong to cases (i) and (ii) of Theorem \ref{cross}. \smallskip

{\it Pisot examples for case (iii). } 
Now we combine two periodic addresses. These examples indicate parameters $s$ where $\nu_s$ does not possess a bounded density, and perhaps is even singular.  We study the peaks of Figure \ref{B7}, and show a close-up of an important part of phase 3 in Figure \ref{B9} below. The curves $y_b(t)$ with $b=3/7=.\overline{011}$ and $b=4/9=.\overline{011100}$ and $y_c(t)$ with $c=8/15=.\overline{1000}$ and $c=16/31=.\overline{10000}$ go through this region. We can determine their four intersection points of type (iii), providing the Pisot parameter $s=.570$ from Table \ref{pigar} and three other Pisot parameters $s=.552, .560, .576.$  Thus we have singularities of the function $\Phi ,$ according to the old result of Erd\"os. 

There is some more information about local dimensions, however, and even for Pisot numbers other than multinacci, such information is far from being obvious \cite{F4,F11,HHM}. The growth rate $\rho$ of a point lying on the intersection of two cycles of length $m, n$  is at least as large as the positive root of $r^{-m}+r^{-n}=1.$ For $3/7$ and $8/15$ we have $m=3,\, n=4$ and $\rho\ge 1.22 .$ The local dimension then is  $\log \frac{2}{\rho}/\log\beta \le .895 ,$ much smaller than all local dimensions for the Fibonacci parameter. For $4/9$ and $8/15$ we obtain the same growth: since $4/9=.\overline{011100}$ has the form $\overline{w(1-w)}$ with $w=011$ we can also take $m=3.$ When we replace $8/15$ by $16/31,$ this means $n=5$ and  $\rho\ge 1.19.$ Local dimensions are between .87 and .89.\smallskip

\begin{figure}[h]
\includegraphics[width=0.999\textwidth]{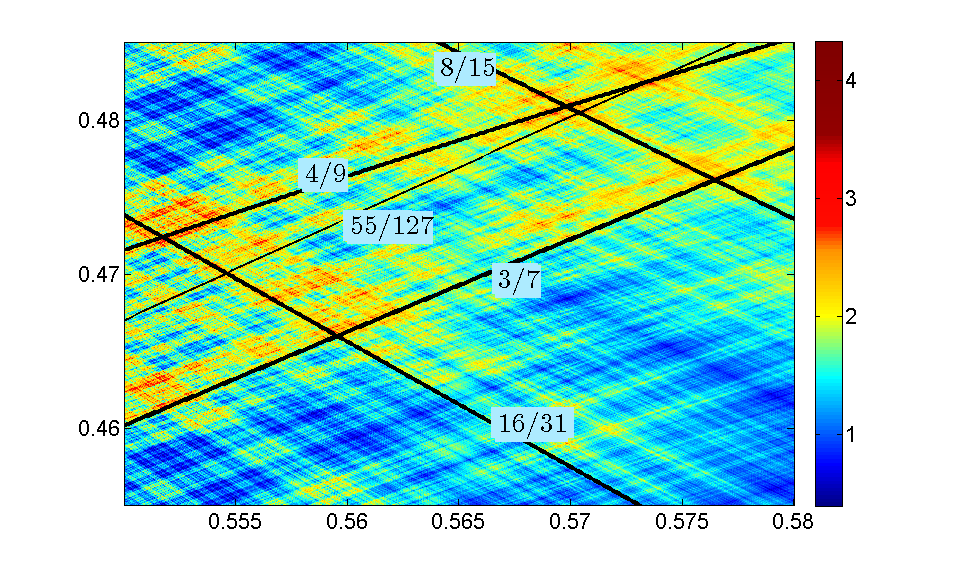}
\caption{$\Phi$ for $0.55\le t\le 0.58,\quad 0.455\le y\le 0.485,$ with address curves of $3/7, 4/9, 8/15, 16/31$ (fat lines) and $55/127$ thin line.  Intersections lead to Pisot and Perron parameters where a density of the Bernoulli measure cannot be bounded.}
\label{B9}
\end{figure}

{\it Are there non-Pisot parameters with singularities?} 
This question was studied by Feng and Wang \cite{FW}, based on work of Peres and Solomyak, and for Salem numbers in greater detail by Feng \cite{F11}.  A main result of \cite{FW} says that a root $t$ of a polynomial of degree $n$ with coefficients $0, \pm 1$ which fulfils $t<2^{-n/(n+1)}$ gives rise to a Bernoulli convolution $\nu_t$ which cannot have a bounded density. Our theorem above allows to find new parameters $t$ and to give an interpretation of the type of singularity which we have: there is a finite orbit of the multivalued map $g_t$ which has growth leading to a local dimension smaller than 1 \cite{BZ}. At the points of this orbit, and at all their preimages under $g_t,$ we must have poles when a density exists. Thus we would have a dense set of poles. These poles can be verified as in Figure \ref{BC1}. 
Conclusions concerning the multifractal spectrum similar to the results in \cite{F4,F11,HHM} can be easily drawn, but  here we confine ourselves to a simple example. \smallskip

{\it A new Perron number with no bounded density. }
Take the address curve of $b=55/127= .\overline{0110111},$ drawn as thin line in Figure \ref{B9}, with formula  $y_b(t)=(t-t^3+t^4-t^7)/(1-t^7).$ Since this $b$ is an itinerary, not a kneading sequence, the curve intersects not only the $y_c(t)$ coming from above, but also some $y_b(t)$ leading upwards. The intersection with the curve of $4/9$ leads to another Pisot parameter at $t=.5735.$  However, the intersection with the curve of $c=16/31$ leads to $t=.5546,$ with minimal polynomial $t^9+t^8+2t^7+t^6+2t^5+t^4+t^3+t-1.$ The number $\beta=1/t$ is Perron, not Pisot. The Feng-Wang result does not apply since $t>2^{-9/10}=.536.$  If there was a polynomial with coefficients $0, \pm 1$ then $2^{-n/(n+1)}$ would be still smaller.

We have the intersection of two periodic orbits with $m=7$ and $n=5,$ and the inequality of the theorem is fulfilled. 
The local dimension can be estimated as above: $\rho\ge 1.1237>1.11\ge 2t$ and  $d_z(\nu_s)= \log \frac{2}{\rho}/\log\beta \le .98 .$  

It should be noted that the intersection of $y_b(t)$ with the curve of $c=8/15$ also yields a Perron parameter, but for this intersection point the inequality in (iii) is not satisfied.
However, there are many kneading sequences $b$ which lead to Perron parameters for which $\nu_t$ has non-trivial multifractal spectrum. This will not be discussed in this introductory note. It seems possible that all poles and local maxima of the function $\Phi$ can be represented as intersections of address curves.

\end{document}